\input amstex

\documentstyle{amsppt}
\magnification=1200
\hcorrection{.25in}
\advance\vsize-.75in

\NoBlackBoxes
\nologo

\def\lap{\triangle}
\def\grad{\nabla}
\def\ip#1#2{\left<#1,#2\right>}
\def\ipmu#1#2{\left<#1,#2\right>_\mu}
\def\iploc#1#2{\left(#1,#2\right)_x}
\def\image{\operatorname{im}}
\def\real{\text{{\bf R}}}

\def\vf{\varphi}

\def\eps{\epsilon}
\def\dom{\Cal D}

\def\Anderson{[{\bf An}]}
\def\Atiyah{[{\bf At}]}
\def\APS{[{\bf APS}]}
\def\Bismut{[{\bf Bi}]}
\def\BottTu{[{\bf BT}]}
\def\BruningLesch{[{\bf BL}]}
\def\Bueler{[{\bf Bu}]}
\def\Chernoff{[{\bf C}]}
\def\CyconEtAl{[{\bf CFKS}]}
\def\Gaffney{[{\bf G}]}
\def\GriffithsHarris{[{\bf GH}]}
\def\Lott{[{\bf L}]}
\def\Melrose{[{\bf M}]}
\def\Prokhorenkov{[{\bf P}]}
\def\ReedSimon{[{\bf RS}]}
\def\Warner{[{\bf Wa}]}
\def\Witten{[{\bf Wi}]}

\topmatter
\title Hodge Theory and Cohomology with Compact Supports \endtitle
\author Ed Bueler \\ Igor Prokhorenkov \endauthor
\leftheadtext{ED BUELER AND IGOR PROKHORENKOV}
\date April 27, 1998 \enddate

\address Ed Bueler --- Mathematical Sciences Research Institute, Berkeley \endaddress

\email ebueler\@msri.org \endemail

\address Igor Prokhorenkov --- Current address:  Kosmonavtov 31/3--94, Lipetsk, 398032, Russia.  After 9/1/98: Department of Mathematics, Texas Christian University, Fort Worth, TX. \endaddress

\email prokhor\@mccme.ru \endemail

\thanks Research at MSRI is supported in part by NSF grant DMS--9701755. \endthanks

\abstract
This paper constructs a Hodge theory of noncompact topologically tame manifolds $M$.  The main result is an isomorphism between the de~Rham cohomology with compact supports of $M$ and the kernel of the Hodge--Witten--Bismut Laplacian $\lap_\mu$ associated to a measure $d\mu$ which has sufficiently rapid growth at infinity on $M$.  This follows from the construction of a space of forms associated to $\lap_\mu$ which satisfy an ``extension by zero'' property.  The ``extension by zero'' property is proved for manifolds with cylindrical ends possessing gaussian growth measures.
\endabstract

\endtopmatter

\document

\heading 1. Introduction \endheading

Let $M$ be a real $C^\infty$ manifold without boundary.  We say $M$ is {\sl topologically tame} if there exists a compact manifold $N$ with boundary such that $M$ is diffeomorphic to the interior of $N$.  For instance, vector bundles over compact base manifolds are topologically tame.  Compact manifolds are trivially topologically tame.  A surface of infinite genus is not topologically tame.

We will construct a Hodge theory of topologically tame noncompact manifolds $M$ in this paper.  First, we choose a metric $g$ such that $(M,g)$ is complete and the ends of $M$  are metric cylinders.  Then we define a particular (Hodge--Witten--Bismut) Laplacian operator $\lap_\mu$ on certain differential forms on $M$.  The subscript $\mu$ denotes an additional choice of measure $d\mu$.  The operator $\lap_\mu$ acts on forms which are square--integrable with respect to the measure $d\mu$.

The essential construct of this paper is that the measure $d\mu$ has certain rapid growth at infinity on $M$ relative to Riemann--Lebesgue measure $dx$: $d\mu = e^{cr^2} dx$, $c>0$.  In this paper we prove (theorem 15):

\proclaim{Main Theorem}  With these choices, $(M^n,g)$, $d\mu$, the $\lap_\mu$--harmonic forms represent the de Rham cohomology with compact supports:
	$$\ker \lap_\mu^p \cong H_{c,\text{{\rm de~R}}}^p(M), \qquad p=0,\dots,n.$$\endproclaim

The skeptical reader may observe that any nonzero $\omega \in \ker\lap_\mu$, as an eigenform of an elliptic operator, {\sl does not} have compact support.  Nevertheless, we will exhibit a (not canonical!)~linear surjection $j: \ker\lap_\mu \to H_{c,\text{{\rm de R}}}(M)$ and a (also not canonical) linear surjection $\Pi: H_{c,\text{{\rm de R}}}(M) \to \ker\lap_\mu$, in the situation covered by the main theorem.

In fact, we will consider a condition on more general measures $d\mu$ and on the domain of their associated operators $\lap_\mu$.  This condition is the extension by zero condition ``property EZ'' defined in the third section, and it is the main technique used in proving theorem 15:

\proclaim{Main Technique}  Choose a ``nice'' relatively compact open submanifold $M_0\subset M$ such that $M_0$ is diffeomorphic to $M$ under $\Psi:M\to M_0$.  Show that if $\omega \in C^\infty(\lap_\mu) = \bigcap_{k\le 1} \dom_{(\lap_\mu)^k}$ then $\left(\Psi^{-1}\right)^*\omega \in \Omega(M_0)$ has smooth extension by zero to $M$.  (``{\rm Property EZ}'')  Show that this condition implies the existence of a surjection $j:\ker \lap_\mu \to H_{\text{\rm c,de~R}}(M)$.\endproclaim

Note that any eigenforms of $\lap_\mu$ are in $C^\infty(\lap_\mu)$, and that they decay very rapidly if $d\mu = e^{c r^2}dx$, $c>0$ as prescribed.  

Property EZ can hold only if $M$ is topologically tame.  

Theorem 13 is that property EZ implies the existence and surjectivity of $j$.

The remainder of the introduction defines $\lap_\mu$ and motivates our choice of $d\mu$.  We also recall the existing Hodge theory of $\lap$, constructed in the usual way from $dx$, in the case of $(M,g)$ with cylindrical ends.  Comparison then shows the role of the densities $e^{\pm r^2}$ from the point of view of the Hodge theory of such $(M,g)$.

Now, the metric $g$ determines an inner product $\iploc\cdot\cdot$ on the tangent space $T_x M$ by definition, which induces an inner product $(\cdot,\cdot)_x$ on $\wedge T_x^* M$, the space of pointwise values of differential forms.  The metric also determines a measure $dx$, by the coordinate formula $dx = \sqrt{\det g_{ij}} dx^1\dots dx^n$.

Assume $M$ is oriented.  It is traditional to define an inner product on the space of smooth ($C^\infty$) compactly--supported differential forms $\Omega_c(M)$ by the formula $\ip\omega\nu = \int_M \iploc\omega\nu dx = \int_M \omega \wedge \star \nu$.  (See \Warner, for instance, for details.)  Here $\star$ is the Hodge duality operator defined in terms of $(\cdot,\cdot)_x$ and the orientation on $M$.

Next we weight the inner product.  Restricting ourselves to smooth, positive densities, if $h\in C^\infty(M)$ then define
	$$d\mu = e^{2h} dx \qquad \text{ and } \qquad \ipmu\omega\nu = \int_M \iploc\omega\nu d\mu.$$

Let $L_\mu^2 \Omega$ be the completion of $\Omega_c(M)$ with respect to $\ipmu\cdot\cdot$.  The exterior derivative $d$ has formal adjoint $\delta_\mu$: $\ipmu{d\omega}\nu = \ipmu\omega{\delta_\mu \nu}$ for $\omega, \nu \in \Omega_c$.  In fact, the formula $\delta_\mu = e^{-2h} \delta e^{2h}$ holds, where $\delta = (-1)^{n(p+1)+1} \star d \star$ acting on $\Omega^p$ is the formal adjoint of $d$ with respect to the usual inner product.

Let $\lap_\mu = d \delta_\mu + \delta_\mu d$.  It is well--known that

\proclaim{1. Theorem}  The operator $\lap_\mu$ defined on $\Omega_c$ is densely--defined and symmetric with respect to $\ipmu\cdot\cdot$, and is a second--order elliptic operator which maps $\Omega_c^p$ to itself.  If $(M,g)$ is complete, we can extend $\lap_\mu$ to a unique nonnegative (self--adjoint) operator in $L_\mu^2 \Omega$ ({\rm \Gaffney, \Chernoff}---see {\rm \Bueler}~for the weighted case).

If $M^n$ is in fact {\rm compact}, then $\lap_\mu$ has compact resolvent (discrete spectrum).  The Hodge decomposition of the smooth forms $\Omega^p(M)$, $p=0,\dots,n$, follows:
	$$\Omega^p(M) = \ker \lap_\mu^p \oplus \image d^{p-1} \oplus \image \delta_\mu^{p+1},\qquad (M \text{ compact})$$
from which it follows that $H_{\text{{\rm de R}}}^p(M) \cong \ker \lap_\mu^p$.\endproclaim

The weighting $d\mu = e^{2h} dx$ on $L^2$ spaces of forms has been considered in connection with the Witten treatment of the Morse inequalities \Witten.  Actually, our $\lap_\mu$ on $L_\mu^2\Omega$ is unitarily--equivalent to the Laplacian usually considered in that context, as explicitly noted in \Bismut.  There is the following unitary equivalence:
$$\matrix L_\mu^2\Omega & @>U>> & L^2\Omega \\
	      \omega            & \mapsto &  e^h \omega \\
		d		&	  & d_h = U d U^{-1} = e^h d e^{-h} \\
  \lap_\mu = d\delta_\mu + \delta_\mu d & & \lap_h = U \lap_\mu U^{-1} = \lap + |dh|^2 + A_h, \endmatrix \tag1$$
where $A_h$ is a certain zero order symmetric operator on $L^2\Omega(M)$ involving the second derivatives of $h$.

We want to consider noncompact topologically tame manifolds $M$.  On such $M$ we can explicitly construct metric and measure which yield nice spectral properties.

\proclaim{2. Theorem}  ({\rm \Bueler}, Theorem 10.2)  Let $M$ be a $C^\infty$, orientable, connected, and topologically tame manifold.  There exists a metric $g$ and for each $c\ne 0$ a measure $d\mu$ with smooth density with respect to $dx$ such that
$$\align (i) \quad &(M,g) \text{ is complete and} \\
	(ii) \quad &\lap_\mu \text{ in } L_\mu^2 \Omega \text{ has compact resolvent (= discrete spectrum).} \endalign$$
Under the chosen metric $g$, the ends of $M$ are cylinders near infinity.  The density of the chosen measure $d\mu$ behaves precisely as $\exp(c r^2)$ near infinity on the ends, for any chosen $c \ne 0$. \endproclaim

A certain Hodge decomposition is one consequence of (ii).  It is sometimes misnamed a {\sl strong} Hodge decomposition.  It says:

\proclaim{3. Proposition}  If $(M^n,g)$ is complete and $\lap_\mu$ has compact resolvent, or merely a spectral gap (i.e. $\sigma(\lap_\mu) \setminus \{0\} \subseteq [\eps,\infty)$ for some $\eps > 0$), then 
	$$L_\mu^2 \Omega^p = \ker \lap_\mu^p \oplus \image d^{p-1} \oplus \image \delta_\mu^{p+1}, \qquad p=0,\dots,n.$$\endproclaim

\demo{Proof}  Classical.  Use the spectral theorem.  See \Bueler~for details.\qed \enddemo

The conclusion of proposition 3 is {\sl not} strong in the sense that there is {\sl a priori} no relation between $\ker\lap_\mu$ and either $H_{\text{{\rm de~R}}}(M)$ or $H_{c,\text{{\rm de~R}}}(M)$ in the noncompact case.  On must produce a separate proof of any correspondence.  We will do so in the case covered by theorem 2---we will show that if $c>0$ then $\ker \lap_\mu \cong H_{c,\text{{\rm de~R}}}(M)$ and if $c<0$ then $\ker \lap_\mu \cong H_{\text{{\rm de~R}}}(M)$.

The widely studied case $d\mu = dx$ provides plenty of examples showing that $\ker\lap$ and $H_{\text{{\rm de~R}}}(M)$ at best stand in a very complicated relation to each other.  The surveys \Anderson~and \Lott~should serve as adequate illustration.

A positive example in the case $d\mu=dx$ is that if $M$ has the additional structure of a Riemannian cover of a compact quotient manifold $M/\Gamma$, then the space $\ker \lap$ is determined by the topology of $M/\Gamma$ and by the group $\Gamma$.  This situation was addressed by Atiyah \Atiyah.  (It should be clear that in general these noncompact $M$, with compact quotient $M/\Gamma$, need not be topologically tame.  Conversely, many topologically tame noncompact $M$ admit no metric under which there exists a group of isometries $\Gamma$ with compact quotient $M/\Gamma$.)

The class of Riemannian manifolds $(M,g)$ we consider in theorems 2 and 15---topologically tame and with cylindrical ends---has been a subject of Hodge--theoretic consideration.  Again take $d\mu=dx$.  The following theorem of Atiyah, Patodi, and Singer applies:

\proclaim{4.  Theorem} ({\rm \APS}, Proposition 4.9)  Assume $(M,g)$ is complete, is topologically tame and has cylindrical ends.  Let $i:\Omega_c(M) \to \Omega(M)$ be the inclusion map.  The space $\ker \lap$ of $L^2$ harmonic forms is isomorphic to the {\rm image of the induced map} $i$ from $H_{c,\text{{\rm de~R}}}(M)$ into $H_{\text{{\rm de~R}}}(M)$. \endproclaim

Although this result establishes a connection between the space of $L^2$ harmonic forms and the (topological) cohomologies, it shows that $\ker\lap$ must vanish in circumstances where there is much cohomology to be represented.  For instance, $\ker\lap^p$ vanishes for all $0\le p \le n$ if $M^n$ is not compact and either $H_{\text{{\rm de~R}}}^p(M)=0$ for $1\le p \le n/2$ or $H_{\text{{\rm de~R}}}^p(M)=0$ for $n/2\le p \le n-1$.  (Recall that $\ker \lap^p \cong \ker \lap^{n-p}$ for any $(M,g)$ and that $H_{\text{{\rm de~R}}}^n(M)=0$ if $M$ is not compact.)

Theorem 4 and the paper \APS~have motivated a family of results on manifolds with cylindrical ends.  R.~Melrose and coworkers (specifically, \Melrose) address the situation of a compact manifold with boundary, and a certain prescribed metric near the boundary which they call an ``exact $b$--metric''.  Section 6.4 of \Melrose~includes a Hodge theorem which may be equivalent to our theorem 15 after mapping the manifold with boundary and exact $b$--metric to a noncompact manifold with cylindrical ends.  However, the technique given in section 3 of this paper applies to manifolds with noncylindrical ends, and furthermore directly addresses the complete noncompact manifold itself.  For instance, the ``boundary condition'' here, that makes $\lap_\mu$ self--adjoint, is exactly the condition that $(M,g)$ is complete.

Now, we are motivated to consider measures with particular decay rates by the following example, based on well--known spectral properties of the quantum harmonic oscillator.

\proclaim{5. Example}  {\rm Let $d\mu^+ = e^{+|x|^2} dx$, $d\mu^- = e^{-|x|^2} dx$ on Euclidean space $\real^n$.  Then $\lap_{\mu^+}$ and $\lap_\mu^-$ have compact resolvent and
	$$\ker\lap_{\mu^-}^p \cong \cases \real, &p=0 \\ 0, &p>0,\endcases \qquad \ker\lap_{\mu^+}^p \cong \cases 0, &p<n \\ \real, &p=n.\endcases$$
Note that
	$$H_{\text{{\rm de~R}}}^p(\real^n) \cong \cases \real, &p=0 \\ 0, &p>0,\endcases \qquad H_{c,\text{{\rm de~R}}}^p(\real^n) \cong \cases 0, &p<n \\ \real, &p=n.\endcases$$}
\endproclaim

In the example above, it is clear that the cases $d\mu^+ = e^{c |x|^2} dx$ for $c>0$ and $d\mu^- = e^{c |x|^2} dx$ for $c<0$ on $\real^n$ differ fundamentally, but we can easily show their close relation as well.  In general, if $d\mu^+ = e^{+2h} dx$ and $d\mu^- = e^{-2h} dx$ on $M$, then from (1)
	$$\lap_{+h} = U_+ \lap_{\mu^+} U_+^{-1} = \lap + |dh|^2 + A_h,$$
	$$\lap_{-h} = U_- \lap_{\mu^-} U_-^{-1} = \lap + |dh|^2 - A_h,$$
are the unitarily equivalent versions of $\lap_{\mu^+}$, $\lap_{\mu^-}$ respectively, both acting on (unweighted) $L^2\Omega$.  Thus the difference between $\lap_{\mu^+}$ and $\lap_{\mu^-}$ can be thought of as determined only by the sign of the potential term $A_h$ (an endomorphism on $\wedge T^* M$), once we push $\lap_{\mu^+}$, $\lap_{\mu^-}$ back into the same space $L^2\Omega$.  On the other hand, note that $L_{\mu^+}^2\Omega$ and $L_{\mu^-}^2\Omega$ contain forms of (in general) very different growth/decay at infinity.

In the rest of this paper, we will drop ``de~R'' from our notation, but it is appropriate to note that our results relating analytic quantities to topological quantities are ``modulo de~Rham's theorem''.  Recall that de~Rham's theorem identifies the de~Rham space $\ker d/\image d$ defined in the $C^\infty$ category with the real cohomology space $H(M;\real)$ defined in the continuous category.

Now, if $\ker \lap_\mu^p$ represents $H^p(M)$ for all $0\le p \le n$ for a given noncompact $M$, then it will not also represent $H_c^p(M)$.  However, as suggested by example 4, measures with inverse densities give a desirable relationship.

\proclaim{6. Proposition}  {\rm Poincar\'e duality.}  If $M^n$ is orientable, then
	$$H^p(M) \cong \left(H_c^{n-p}(M)\right)^*.$$
{\rm Weighted Poincar\'e duality.}  If $(M^n,g)$ is orientable, and $d\mu^+ = e^{+2h} dx$, $d\mu^- = e^{-2h} dx$, then
	$$\ker\lap_{\mu^-}^p \cong \ker\lap_{\mu^+}^{n-p}.$$\endproclaim

\demo{Proof}  For Poincar\'e duality, see \BottTu.  For this weighted version, see \Bueler. \qed \enddemo

It follows that if we find a measure $d\mu = e^{2h} dx$ for which $\ker\lap_{\mu}^p \cong H_c^p(M)$, then defining $d\mu^- = e^{-2h} dx$, we get (for free!) $\ker\lap_{\mu^-}^p \cong H^p(M)$.  Theorem 15 shows that this {\sl is} the situation that holds for a topologically tame manifold, having chosen $g$ and $d\mu$ for $c>0$ and $c<0$, respectively, as in theorem 2.

The methods of this paper are closely related to those in \Prokhorenkov, where, in particular, the isomorphism between the $L^2$ cohomology of a vector bundle $E$, corresponding to a density $e^{c|x|^2}$ for $c>0$, and the compact cohomology $H_c(E)$ was proven, and then used in a proof of the Morse--Bott inequalities.

However, the particular form of definitions 7, 9 and 11 and of theorem 13 are motivated by the following situation.  If $(M,g)$ is complete and has Ricci curvature bounded below, then the usual Laplace--Beltrami operator $-{1\over 2} \lap$ on $M$ has a unique heat kernel $\rho$.  On Euclidean space $\real^n$, $\rho_t(x,y) = (2\pi t)^{-n/2} \exp(-|x-y|^2/2t)$.  On $M$ with Ricci bounded below, if $t>0$ and $x_0\in M$ then $d\mu^- \equiv \rho_t(x_0,x) dx$ is a probability measure on $M$.  It is conjectured in \Bueler~ that $\ker \lap_{\mu^-} \cong H^p(M)$ for this heat kernel measure $d\mu^-$, and certain evidence is given there.  It is hoped that the methods of this paper will apply to $d\mu^+ \equiv \left(\rho_t(x_0,x)\right)^{-1} dx$, in which case the conjecture would hold true.  The conjecture implies an existence conjecture: if $b_p>0$ then $\lap_{\mu^-}^p$ has nonvanishing kernel.

We thank Rostislav Matveyev for conversations and help with lemma 12.  The second author thanks Maxim Braverman for valuable disscussions about the subject of section 5.

\heading 2. Forms of Very Rapid Decay $C_\mu^\infty$ \endheading

\proclaim{7. Definition}  Assume $(M,g)$ is complete, so $\lap_\mu$ is self--adjoint.  Define
	$$C_\mu^\infty \equiv C^\infty(\lap_\mu) = \bigcap_{k\ge 1} \dom_{(\lap_\mu)^k}.$$ \endproclaim

Clearly $\Omega_c(M) \subseteq C_\mu^\infty$, so $C_\mu^\infty$ is never empty.  If $\lap_\mu$ has eigenforms, they too are in $C_\mu^\infty$: $\lap_\mu\omega=\lambda\omega$ clearly implies $\omega\in \bigcap_k \Cal D_{(\lap_\mu)^k}$.

The ellipticity of $\lap_\mu$ implies that $C_\mu^\infty \subseteq \Omega(M)$, i.e. that $C_\mu^\infty$ forms are smooth.  If $M$ is compact, then $C_\mu^\infty = \Omega(M)$.  If $M$ is not compact, then $C_\mu^\infty$ includes constraints on growth at infinity, and in fact $\Omega_c(M) \subset C_\mu^\infty \subset \Omega(M)$ are both proper inclusions.

It is appropriate to think of $C_\mu^\infty$ as some kind of generalization to forms of the Schwartz space of functions, with the growth bounds determined by $d\mu$.  More accurately, $e^h C_\mu^\infty = C^\infty(\lap_h)$ generalizes Schwartz space---see section V.

Growth constraints are included in definition 7.  Consider the example:

\proclaim{8. Example} {\rm Let $M=\real^n$ with the Euclidean metric, and let $d\mu=e^{+x^2} dx$.  The unitarily equivalent operator $\lap_h$ on $L^2\Omega$, corresponding to $\lap_\mu$ on $L_\mu^2\Omega$, in this case is the harmonic oscillator Hamiltonian $\lap_h^p = \lap^p + |x|^2 + n -2p$ (on $L^2\Omega^p$).

Considering only functions, it is well known that $C^\infty(\lap_h^0)$ corresponds to the Schwarz space $\Cal S(\real^n)$.  In fact, the family of seminorms $\{\|f\|_k = \|(\lap_h^0)^k f\|_{L^2}\}$ is equivalent to the usual family of seminorms $\{\|f\|_{\alpha,\beta} = \|x^\alpha D^\beta f\|_{L^2}\}$ ({\rm \ReedSimon}, Appendix to V.3).  If follows that $f\in C^\infty(\lap_h^0)$ iff $\|f\|_k <\infty$ for all $k$ iff $f\in \Cal S(\real^n)$.  Therefore $g \in C_\mu^\infty$ iff $g = e^{-|x|^2/2} f$ for $f\in \Cal S(\real^n)$.   This shows that functions in $C_\mu^\infty$ have ``very rapid'' decay.}\endproclaim

\proclaim{9. Definition}   $H_\mu$ is the cohomology of $(C_\mu^\infty,d)$:
	$$H_\mu^p \equiv {\ker d|_{C_\mu^{\infty,p}} \over \image d|_{C_\mu^{\infty,p-1}}}, \qquad p=0,\dots,n.$$ \endproclaim

The closure of $\image d|_{C_\mu^\infty}$ is {\sl not} taken in our definition of $H_\mu$.  Therefore $H_\mu \cong \ker \lap_\mu$ if $\lap_\mu$ has a spectral gap:

\proclaim{10. Corollary}  The Hodge decomposition of proposition 3 can be smoothed:  If $\lap_\mu$ has a spectral gap, then
	$$C_\mu^\infty = \ker \lap_\mu \oplus \image d|_{C_\mu^\infty} \oplus \image \delta_\mu|_{C_\mu^\infty}$$
so $\ker \lap_\mu \cong H_\mu.$\endproclaim

\demo{Proof}  First, clearly $\ker\lap_\mu \subset C_\mu^\infty$.  Thus from
proposition 3, 
	$$C_\mu^\infty = \ker \lap_\mu \oplus C_\mu^\infty \cap \left(\image d \oplus \image \delta_\mu\right).$$

Suppose $\omega \in C_\mu^\infty$, $\omega = d\phi + \delta_\mu\psi$, with $\phi$,$\psi$ not necessarily smooth.  Recall that if $\nu\in\ker\lap_\mu$ then $d\nu = \delta_\mu\nu = 0$.  We have
	$$\ipmu\nu\omega = \ipmu\nu{d\phi} + \ipmu\nu{\delta_\mu\psi} = \ipmu{\delta_\mu\nu}\phi + \ipmu{d\nu}\psi = 0$$
for all $\nu \in \ker\lap_\mu$.  So $\omega \in (\ker\lap_\mu)^\perp$.  Since $\lap_\mu$ has a spectral gap, $(\lap_\mu)^{-1}$ is a bounded operator on $(\ker\lap_\mu)^\perp$.  Thus $\omega = \lap_\mu \lap_\mu^{-1} \omega = d \left(\delta_\mu \lap_\mu^{-1} \omega\right) + \delta_\mu \left(d \lap_\mu^{-1} \omega\right)$.  That is, we can take $\phi=\delta_\mu \lap_\mu^{-1} \omega$, $\psi=d\lap_\mu^{-1}\omega$ in $C_\mu^\infty$.\qed \enddemo

\heading 3. Extension--by--Zero and Compactly--Supported Cohomology \endheading

We start with an obvious construction on a topologically tame manifold.

Suppose $M$ is topologically tame and is equipped with a metric $g$ making $(M,g)$ a complete metric space.  Then we can think of $M$ as including a finite diameter part carrying all the interesting topology, and a finite list of infinite--length ends which are topological products.  In fact, let $\Phi : \operatorname{int} N \to M$ be a diffeomorphism.  Let $\hat U \subset N$ be any open neighborhood of $\partial N$ which is small enough so that $\hat U \cong \partial N \times (0,\infty]$.  The closure of $\hat U$ in $N$ is compact.  Decompose $\hat U$ into finitely--many connected pieces, so $\hat U = \bigsqcup_i Q_i \times (0,\infty]$, with each $Q_i$ compact without boundary.  The boundary $\partial N$ is identified with points $(q\in Q_i, r=\infty)$ in this representation of $\hat U$.   Define $U = \bigsqcup_i Q_i \times (0,\infty)$.  Let $M_0 = M\setminus \overline{\Phi(U)}$, so $M_0$ has compact closure in $M$ and finite diameter, and let $E_i = \Phi(Q_i \times (0,\infty))$.  We can decompose $M = \overline{M_0} \sqcup \left(\bigsqcup_i E_i\right)$.  Note that $M_0$ is diffeomorphic to $M$ itself.

Now we consider the induced mapping $(\Psi^{-1})^*$, for $\Psi : M \to M_0$ a diffeomorphism, taking ``very--rapid--decay'' forms $C_\mu^\infty$ back to $M_0$.

\proclaim{11. Definition}  {\rm Property EZ.}  Let $M$ be topologically tame, assume $(M,g)$ is complete, and suppose a measure $d\mu=e^{2h} dx$ is given.  Consider $M_0 \subset M$ an open submanifold of $M$ with compact closure, such that $\partial M_0 \subset M$ is a smooth closed submanifold, and such that there exists a diffeomorphism $\Psi:M \to M_0$.  Such $M_0$ exist since $M$ is topologically tame.  We say $d\mu$ has {\rm Property EZ} if there exists $M_0$ and $\Psi$ so that every $\omega \in (\Psi^{-1})^* C_\mu^\infty$ has a smooth extension by zero to all of $M$. \endproclaim

From now on, we will write $H_c(M)$ for the cohomology of the complex $(\Omega_c(M),d)$, i.e.~for the de~Rham cohomology with compact supports.  Note also that the cohomology spaces $H_c(M_0)$, $H_c(M)$ in the following lemma are finite--dimensional because $M$ is topologically tame.

\proclaim{12. Lemma}  Let $M$ be a topologically tame manifold.  Suppose $M_0\subset M$ is an open submanifold with compact closure, $\Psi:M\to M_0$ a diffeomorphism.  Assume $\partial M_0$ is a smooth (closed) submanifold of $M$.  Let $i:\Omega_c(M_0) \to \Omega_c(M)$ be the map which extends forms by zero.  Then $i$ induces an isomorphism $i:H_c(M_0) \to H_c(M)$.\endproclaim

\demo{Proof}  We construct a family $\Psi_\eps$, $0< \eps < 1$, of $M \to M_0$ diffeomorphisms which fix ``large portions'' of $M_0$.  First, note that $\Psi(M_0)$ is diffeomorphic to $M$, since $\Psi|_{M_0} \circ \Psi$ is a diffeomorphism.  Thus it is clear that $M\setminus \Psi(M_0)$ is diffeomorphic to $\partial M_0$ times an interval.  Let $\alpha : \partial M_0 \times [-1,\infty) \to M\setminus\Psi(M_0)$ be such a diffeomorphism, chosen so that $\alpha(\partial M_0\times [-1,0)) = M_0 \setminus \Psi(M_0)$ and $\alpha(\partial M_0 \times [0,\infty)) = M\setminus M_0$, so $\alpha(\partial M_0\times \{0\}) = \partial M_0$.  Let $M_\eps = M_0 \setminus \alpha(\partial M_0\times (-\eps,\infty))$ for $0<\eps<1$, and also define $\beta_\eps:[-1,\infty) \to [-1,0)$ so that $\beta_\eps(\xi) = \xi$ for $-1\le \xi\le -\eps$, and $\beta_\eps$ is smooth and monotone increasing with $\lim_{\xi\to\infty} \beta_\eps(\xi) = 0$.

Define $\Psi_\eps(x) = x$ for $x\in \Psi(M_0)$.  For $x\in M\setminus\Psi(M_0)$, let $(y,\xi) = \alpha^{-1}(x)$, $y\in\partial M_0$, $\xi\in [-1,\infty)$.  Define $\Psi_\eps(x) = \alpha(y,\beta_\eps(\xi))$.  Then $\Psi_\eps$ fixes $M_\eps$ and is a diffeomorphism $M\to M_0$ for each $0<\eps<1$.

Now, $\Psi_\eps$ induces $\Psi_\eps^*:\Omega_c(M_0)\to \Omega_c(M)$ on forms and an isomorphism $\Psi_\eps^*:H_c(M_0) \to H_c(M)$.

If $\alpha\in H_c(M_0)$, $\alpha\ne 0$, we will show $i\alpha\ne 0$ in $H_c(M)$.  This suffices, since $H_c(M_0)$, $H_c(M)$ are finite--dimensional.  Choose any representative $\omega$, so $[\omega]=\alpha$.  Since $\omega$ has compact support in $M_0$, there exists $\eps>0$ such that $\operatorname{supp} \omega \subset M_\eps$.  Then $\Psi_\eps^* \omega = i \omega$.  But $i\alpha=[i\omega]=[\Psi_\eps^*\omega]=\Psi_\eps^*\alpha\ne 0$ in $H_c(M)$ since $\Psi_\eps^*$ is an isomorphism.\qed \enddemo

\proclaim{13. Theorem}  Suppose $M^n$ is topologically tame, assume $(M,g)$ is complete, and suppose the measure $d\mu = e^{2h} dx$ is given.  If $d\mu$ has property EZ, then there exists a map $j:C_\mu^\infty \to \Omega_c(M)$, which restricts to a surjection
	$$j:H_\mu \to H_c(M).$$
\endproclaim

\demo{Proof}  Let $M_0$, $\Psi$ be as in the definition of ``property EZ''.  Define the map $j$ to take $\omega \in C_\mu^\infty$ to the extension by zero of $(\Psi^{-1})^*\omega$.  Let $k:\Omega_c(M_0) \to C_\mu^\infty$ be extension by zero.  Then $i$ factors through $j$: $i=j \circ k$.

By the lemma, $i$ is an isomorphism in cohomology.  At the level of cohomology,
	$$H_c(M_0) @>k>> H_\mu @>j>> H_c(M), \qquad i = j \circ k \text{ isomorphism},$$
it follows that $k$ is an injection and $j$ is a surjection. \qed \enddemo

Recalling Poincar\'e duality, we see that if $d\mu$ on $(M^n,g)$ has property EZ then
	$$\dim H_\mu^p \ge \dim H_c^p(M) = \dim H^{n-p}(M) \equiv b_{n-p}.$$

We do not exclude the possibility that $\dim H_\mu^p = \infty$, nor do we assume that $H_\mu^p$ has a Hilbert space structure.  Also, it is not assumed that $\lap_\mu$ has a spectral gap.

There is an available map going the other direction from $j$.  Consider the orthogonal projection $\Pi : L_\mu^2\Omega \to \ker \lap_\mu$.  This map clearly restricts to $\Pi: \Omega_c(M) \to \ker\lap_\mu$.  We also use $\Pi$ to denote the map on compactly--supported cohomology:

\proclaim{14. Definition}  Assume $(M,g)$ is complete.  Define
$$\Pi : H_c(M) \to \ker \lap_\mu \qquad \text{ by } \qquad [\omega] \mapsto \Pi(\omega).$$\endproclaim

Recall that if $M$ is complete, then $\ipmu{d\nu}\varphi = \ipmu\nu{\delta_\mu\varphi} = 0$ for all $\nu \in \Omega_c(M)$ and $\varphi\in\ker\lap_\mu$.  This shows that $\Pi$ is well--defined.

If $d\mu$ has property EZ, then the following method may also apply.  First, the natural map $\eta:\ker\lap_\mu \to H_\mu$, given by $\omega\mapsto [\omega]$ is always injective if $M$ is complete, since if $\omega=d\nu$ for $\nu\in C_\mu^\infty$ then $\ipmu\omega\varphi = \ipmu\nu{\delta_\mu\varphi} = 0$ for all $\varphi\in \ker\lap_\mu$.  If $d\mu$ has property EZ, and if the composition
	$$\ker \lap_\mu @>\eta>> H_\mu @>j>> H_c(M) @>\Pi>> \ker \lap_\mu$$
is surjective, then $\Pi$ is surjective.  This reduces to checking whether an explicit endomorphism $\Pi\circ j\circ\eta$ on $\ker \lap_\mu$, depending only on $\Psi$ and $d\mu$, is surjective.

Whether $\Pi$ is surjective or not is {\sl a priori} independent of the existence of a spectral gap for $\lap_\mu$.

In section 4 we will see that both $j$ and $\Pi$ are surjective in the case described by theorem 2 if $c>0$, and thus that $H_\mu \cong H_c(M)$.  That is, we prove that if: (i) $M$ is topologically tame, (ii) the ends of $(M,g)$ are cylinders, and (iii) $d\mu$ is a measure with density $\exp(cr^2)$, $c>0$ on the ends, then $d\mu$ has property EZ.  It follows that $j$ is surjective, and we prove $\Pi$ is surjective by the method described above.

\heading 4.  Cylindrical Ends, Gaussian Growth Measure \endheading

Let us return to the setting of theorem 2.  Let $M$ be topologically tame, and define $g$ as a metric under which $(M,g)$ is complete and the ends of $M$ are cylinders.  Define $d\mu$ as in theorem 2 as well, with $c>0$: $d\mu=e^{cr^2}dx$.  (We will more explicitly define $g$ and $d\mu$ in a moment.)  In this section we prove that $d\mu$ on $M$ has property EZ, and in fact we prove that both $j$ and $\Pi$ of the previous section are surjections.  As we will see, the proof of property EZ ultimately depends on a fundamental estimate for $C_\mu^\infty$ (actually $C_h^\infty$) forms on $(M,g)$.  The estimate is found in theorem 16 and corollary 17 of section 5.

\proclaim{15. Theorem}  Let $c>0$.  Let $M$, $g$, $d\mu = e^{cr^2} dx$ be as in theorem 2.  Then in addition to (i) and (ii), there exist surjections $\Pi$ and $j$ such that
	$$H_c(M) @>\Pi>> \ker\lap_\mu \cong H_\mu @>j>> H_c(M).$$
It follows that 
	$$H_\mu \cong H_c(M).$$\endproclaim

\demo{Proof that $j$ is surjective}  We show that $d\mu$ has property EZ, so that by theorem 13, $j:H_\mu \to H_c(M)$ is a surjection.

Since $M$ is topologically tame, $M  = M_0 \cup (\bigsqcup_i E_i)$, with $M_0$ relatively compact, $E_i = (0,\infty) \times Q_i$, and $Q_i$ compact without boundary.  Let $g$ be a product metric on each $E_i$ and otherwise choose $g$ arbitrarily.

Fix an end $E = (0,\infty) \times Q$, and let $Q=\bigcup_j U_j$ be a finite cover of $Q$ by coordinate neighborhoods.  Let $(x^1, x^2, \dots, x^n) = (x^1, \bar x) = (r, \bar x)$ be coordinates on $(0,\infty) \times U_j$.

Define $d\mu=e^{cr^2}dx = e^{c(x^1)^2}dx$ on $(0,\infty) \times U_j$, therefore on $E$.  Extend $d\mu$ to all of $M$ in the obvious way, so $d\mu = e^{+2h}dx$ on $M$, where $h\in C^\infty(M)$ and $h= cr^2$ on each end, for the given $c>0$.

Define $\Psi$ on $E$ as follows.  Let $M_1 = M_0 \cup (\bigsqcup_i (0,1)\times Q_i)$, with coordinates $(s,\bar y)$ on $(0,1)\times U_j$.  We define $\Psi^{-1}:M_1\to M$ by defining $\Psi^{-1} = \operatorname{id}$ on $M_0$ and by defining $\Psi^{-1}:(0,1)\times U_j \to (0,1)\times U_j$ as follows: $\Psi^{-1}(s,{\bar y}) = ((1-s)^{-1},{\bar y})$ if $1/2 \le s < 1$, and smoothly and monotonically interpolate $\Psi^{-1}$ for $0\le s\le 1/2$.  Clearly $\Psi^{-1}$ is a diffeomorphism and $M_1 = \Psi(M)$ has compact closure.

To prove that $d\mu$ has property EZ, it suffices to prove that if $\omega\in C_\mu^\infty$ and $\operatorname{supp} \omega \subset E$, then $\left(\Psi^{-1}\right)^* \omega$ has smooth extension by zero to $E$.  That is, if $\vf_0$,$\{\vf_i\}$ is a partition of unity subordinate to $M  = M_1 \cup (\bigsqcup_i E_i)$ then $\omega \in C_\mu^\infty$ iff $\vf_i\omega \in C_\mu^\infty$ for each $i$.

Now we arrange to use the properties of Schwartz forms on an end $E$---see section V---to prove the extension by zero on $E$ for $\left(\Psi^{-1}\right)^* \omega$.

In fact, recall that $\lap_h = e^h \lap_\mu e^{-h} = \lap + |dh|^2 + A_h$ on $L^2 \Omega$ is unitarily--equivalent to $\lap_\mu$ on $L_\mu^2\Omega$.  Here $h={c\over 2} r^2$.  Let $\nu \in C_h^\infty = \bigcap_{k\ge 1} \dom_{(\lap_h)^k}$.  Then $e^{-h} \nu \in C_\mu^\infty$.  Suppose $\nu = \sum_{1\notin J} u_J dr \wedge d{\bar x}^J + \sum_{1\notin I} u_I d{\bar x}^I$ on $(0,\infty)\times U_j$.  Then since ${\partial r \over \partial s} = (1-s)^{-2}$,
$$\align \omega' \equiv \left(\Psi^{-1}\right)^* \left(e^{-h} \nu\right) = &e^{-{c\over 2} (1-s)^{-2}} (1-s)^{-2} u_J((1-s)^{-1},{\bar y}) ds\wedge d{\bar y}^J \\
	&\quad + e^{-{c\over 2} (1-s)^{-2}} u_I((1-s)^{-1},{\bar y}) d{\bar y}^I\endalign$$
is in $\Omega(M_1)$.

\def\ppx#1{{\partial \over \partial x^{#1}}}

Denote as $\ppx j : \Omega^p((0,\infty)\times U_j) \to \Omega^p((0,\infty)\times U_j)$ the operator defined by $\ppx j \omega = \sum {\partial \omega_I \over \partial x^j} dx^I$ if $\omega = \sum \omega_I dx^I$.  That is, $\ppx j$ differentiates the coefficients of $\omega$ with respect to $x^j$.  If $\alpha = (\alpha_1,\dots,\alpha_n) \in \text{\bf Z}_+^n$ is a multiindex, define $\ppx{}^\alpha = \ppx 1^{\alpha_1} \dots \ppx n^{\alpha_n}$.

Clearly, if 
	$$|r^k \left(\partial / \partial x\right)^\alpha \nu(r,\bar x)| \to 0 \text{ \quad as } r \to \infty \tag2$$
for every $k\ge 0$ and every multiindex $\alpha$, then it follows that
	$$|\left(\partial / \partial y\right)^\alpha \omega'(s,\bar y)| \to 0 \text{ \quad as } s \to 1^- \tag3$$
for every $\alpha$, where $(y^1,y^2,\dots,y^n) = (y^1,\bar y) = (s,\bar y)$ are coordinates on $(0,1) \times U_j$.

But (2) is exactly the conclusion of corollary 17, and (3) is exactly what is needed to prove property EZ on the end $E$.  It follows from theorem 13 that $j: H_\mu \to H_c(M)$ is surjective.

({\smc Remark}:  We have proven somewhat more rapid decay in (3) than the minimum needed to prove property EZ.  However, the degree to which the measure $d\mu$ can be changed to a measure $d\mu'$ with slower growth rate $e^{2h'}$ depends on both the growth rate $e^{2h'}$ itself and the corresponding potential $V = |dh'|^2 + A_{h'}$ which plays the essential role in section 5.)\enddemo

\def\wb{\overline{\omega}}
\demo{Proof that $\Pi$ is surjective}  Let $\{\omega^1,\dots,\omega^N\}$ be an orthonormal basis for $\ker \lap_\mu$.  Note $N< \infty$ since $\lap_\mu$ has compact resolvent.  For each $\omega^i$, we will find $\nu \in \Omega_c(M)$ such that $d\nu=0$ and $\Pi(\nu)=\omega^i$.

Choose $\eps=\eps(N)>0$ so that any $N\times N$ real matrix $(a_{ij})$ with $|a_{ij}-1| < \eps$ for all $i$ and $|a_{ij}|<\eps$ for $i\ne j$ is invertible.  For instance, $\eps = N^{-1}$ suffices.

In what follows we will make a finite number of choices for $R$, using the finite--dimensionality of $\ker\lap_\mu$ and the finite list of ends $E_i$, and we assume that the maximum of these $R$ is used as needed.

For any $R\ge 0$, define $\Psi_R$ analogously to $\Psi$ on a particular end $E$: $\Psi_R^{-1}$ is the identity on $\{s\le R\} = M_R$, $\Psi_R^{-1}(s,{\bar x}) = (R+ (1-(s-R))^{-1},{\bar x})$ for $R+1/2 \le s < R+1$, and smoothly and monotonically interpolate $\Psi^{-1}$ for $R\le s\le R+1/2$.  Thus $\Psi_R$ is a diffeomorphism and $\Psi_R(M) = M_{R+1}$.  We may assume $|d\Psi_R^{-1}| \le 4$ on $R\le s\le R+1/2$.  

Next, we use property EZ: define $j_R$ on forms $\vf \in C_\mu^\infty$ as the extension by zero of $(\Psi_R^{-1})^*\vf$.  Let $\wb^i = j_R \omega^i$.  Each $\wb^i \in \Omega_c(M)$ and $d\wb^i = (\Psi_R^{-1})^* d\omega^i = 0$, so $[\wb^i] \in H_c(M)$.

Clearly $|\ipmu{\wb^i}{\omega^j}| \le \|\wb^i-\omega^i\|_\mu$ if $i\ne j$ and $|\ipmu{\wb^i}{\omega^i} - 1| \le \|\wb^i-\omega^i\|_\mu$.  We will show that we can choose $R\ge 0$ so that $\|\wb^i-\omega^i\|_\mu < \eps$.  It will follow that $A = \left(\ipmu{\wb^i}{\omega^j}\right)$ is invertible.

Note $\wb^i|_{M_R} = \omega^i|_{M_R}$. Thus
	$$\|\wb^i-\omega^i\|_\mu^2 = \int_{R<s} |\wb^i-\omega^i|^2 d\mu \le 2 \int_{R<s<R+1} |\wb^i|^2 d\mu + 2 \int_{R<s} |\omega^i|^2 d\mu.$$
Now,
$$\align  \int_{R<s<R+1/2} |\wb^i|^2 d\mu &\le \int_{R<s<R+1/2} \left|{\partial r\over\partial s}\right|^2 |\omega^i(r(s),{\bar y})|^2 e^{cs^2} ds d{\bar y} \\
		&\le 4^2 \int_{R<r} |\omega^i(r,{\bar x})|^2 e^{cr^2} dr d{\bar x}\endalign$$
so the essential estimate is:
$$\align  \int_{R+1/2<s<R+1} |\wb^i|^2 d\mu &= \int_{R+1/2<s<R+1} \left(1-(s-R)\right)^{-2} |\omega^i(r(s),{\bar x})|^2 e^{cs^2} ds d{\bar x} \\
		&= \int_{R+2<r<\infty} |\omega^i(r,{\bar x})|^2 \exp\left(c(R+1 - (r-R)^{-1})^2\right) dr d{\bar x} \\
		&\le \int_{R+2<r<\infty} |\omega^i(r,{\bar x})|^2 \exp\left(c(R+1)^2\right) dr d{\bar x} \\
		&\le \int_{R+1<r} |\omega^i(r,{\bar x})|^2 e^{cr^2} dr d{\bar x}.\endalign$$
It follows that
	$$\|\wb^i-\omega^i\|_\mu^2 \le 36 \int_{R<r} |\omega^i|^2 d\mu.$$
Choose $R\ge0$ so that $\|\wb^i-\omega^i\|_\mu < \eps = \eps(N)$ for every $i=1,\dots,N$.  Then $A = \left(\ipmu{\wb^i}{\omega^j}\right)$ is invertible.

But
	$$\Pi(\wb^i) = \sum_j \ipmu{\wb^i}{\omega^j} \omega^j = A_{ij} \omega^j,$$
so if $\nu=A^{-1}_{ij} \wb^j$ then
	$$\Pi(\nu) = A^{-1}_{ij} \Pi(\wb^j) = A^{-1}_{ij} A_{jk} \omega^k = \omega^i.$$
We have shown that $\Pi$ is surjective, as in the method described in section 3.\qed \enddemo

\heading 5. Schwartz Forms on a Cylindrical End \endheading

Let $E = (0,\infty) \times Q$, where $Q$ is compact without boundary.  Give $E$ a product metric, with $(0,\infty)$ given the Euclidean metric.  (Then $E$ is not complete, but we will compute here as though it is, since for our purposes $E$ is an end of a complete manifold.)  Let $Q=\bigcup_j U_j$ be a finite cover of $Q$ by coordinate neighborhoods.  Let $(x^1, x^2, \dots, x^n) = (x^1, \bar x)$ be coordinates on $(0,\infty) \times U_j$.

Let $h={c\over 2} (x^1)^2$ for $c>0$.

We define $d_h = e^h d e^{-h}$, $d_h^* = e^{-h} d e^h$, and 
	$$\lap_h = d_h d_h^* + d_h^* d_h = \lap + |dh|^2 + A_h,$$
all acting in $L^2 \Omega$ (with measure $dx = \sqrt{\det g_{ij}} dx^i dx^j$ as usual).  Then $\lap_h$ on $L^2 \Omega$ and $\lap_\mu$ on $L_\mu^2\Omega$ are unitarily--equivalent (\Bismut).

Denote $C_h^\infty = \bigcap_{k\ge 1} \dom_{(\lap_h)^k}$.  Clearly, $\omega \in C_h^\infty$ iff $e^{-h} \omega\in C_\mu^\infty$.  We will show that forms in $C_h^\infty$ decay as fast as Schwartz functions on a cylindrical end $E$ as $x^1\to\infty$.

Define $\partial_j : \Omega^p((0,\infty)\times U_j) \to \Omega^p((0,\infty)\times U_j)$ by $\partial_j \omega = \sum {\partial \omega_I \over \partial x^j} dx^I$ if $\omega = \sum \omega_I dx^I$.  That is, $\partial_j$ differentiates the coefficients of $\omega$ with respect to $x^j$.  If $\alpha = (\alpha_1,\dots,\alpha_n) \in \text{\bf Z}_+^n$ is a multiindex, let $\partial^\alpha = \partial_1^{\alpha_1} \dots \partial_n^{\alpha_n}$.

We prove the following result on forms in $C_h^\infty$:

\def\bbN{\text{\bf N}}
\proclaim{16. Theorem}  If $\omega \in C_h^\infty$ is supported in $E$, then for any $k \in \bbN$ and any multiindex $\alpha$,
	$$\|(x^1)^k \partial^\alpha \omega\|_{L^2\Omega} < \infty.\tag4$$
\endproclaim

Before proving (4), we note that it implies a pointwise result as well.  That is, Sobolev inequalites (or even explicit Schwartz space calculations as in \ReedSimon, Appendix to V.3) imply:

\proclaim{17. Corollary}  If $\omega \in C_h^\infty$ is supported in $E$, then for any $k \in \bbN$, any multiindex $\alpha$, and any $\bar x_0 \in U_j$,
	$$|(x^1)^k \partial^\alpha \omega|(x^1,\bar x_0) \to 0 \text{ \quad as } x^1 \to \infty.\tag5$$
\endproclaim

\noindent {\bf 18. Remark.}  Properties (4) and (5) do not depend on the choice of local coordinates $\bar x$ on $Q$, since all transition functions and their derivatives are bounded on compact $Q$.

Also, since $Q$ is compact, by using a partition of unity on $Q$ it is enough to prove theorem 16 only for forms which have support completely inside one of the coordinate neighborhoods $(0,\infty)\times U_j$.  Moreover, since inequality (4) is trivially true for smooth forms with compact support, we can assume (without loss of generality) that $\omega\equiv 0$ for $x^1\le 1$.

\medskip
The proof of theorem 16 will proceed in several steps, starting with following lemmata.

\proclaim{19. Lemma}  {\rm Weitzenbock's formula for $\lap_h$.}  Suppose $\omega$ is a smooth form on $E$, then
	$$\lap_h \omega = \left(\sum_{i=1}^n \grad_i^* \grad_i + R(\bar x) + c^2(x^1)^2 \pm c\right)\omega,\tag 6$$
where $R(\bar x):\wedge^p T_{\bar x}^* E \to \wedge^p T_{\bar x}^* E$ is a pointwise curvature  endomorphism, and $\grad_i = \grad_{\partial/\partial x^i}$.\endproclaim

\demo{Proof}  See \CyconEtAl, Chapter 12.4. \qed \enddemo

\proclaim{20. Lemma}  {\rm G\"arding's inequality.}   Let $\omega$ be a smooth form, suppose $\operatorname{supp}\omega \subset (1,\infty)\times U_j$ and suppose $\omega$, $\lap_h\omega\in L^2\Omega((1,\infty)\times Q)$.  Then
	$$\sum_{i=1}^n \|\grad_i \omega\|^2 + \|x^1 \omega\|^2 \le C\left(\ip{\lap_h \omega}{\omega }+ |\omega\|^2 \right).\tag7$$\endproclaim

\demo{Proof}  Let $J(\tau):(0,\infty)\to \real_+$ be a smooth non--increasing function such that $J(\tau)\equiv 1$ for $\tau\le 2$ and $J(\tau)\equiv 0$ for $\tau \ge 3$.  We define a family of cut--off functions
	$$J_t:(0,\infty)\to \real_+$$
by
	$$J_t(x^1,\bar x) = J\left({x^1\over t}\right).$$

We multiply both sides of (6) by $J_t^2\omega$ and integrate by parts:
$$\ip{\lap_h \omega}{J_t^2 \omega} =\sum_{i=1}^n \ip{\grad_i\omega}{\grad_i(J_t^2\omega)} + \ip{R(\bar x)\omega}{J_t^2\omega} + c^2 \|x^1 J_t \omega\|^2 \pm c\|J_t\omega\|^2.\tag8$$
For each $i$ we have
\def\ips#1#2{\big<#1,#2\big>}
	$$\ip{\grad_i\omega}{\grad_i(J_t^2\omega)} = \|J_t\grad_i\omega\|^2 + 2\ips{J_t {\partial J_t\over \partial x^i}\grad_i \omega}{\omega}.$$
Therefore, substituting back into (8) we get
$$\align \sum_{i=1}^n &\|J_t\grad_i\omega\|^2 + c^2\|x^1 J_t \omega\|^2 \tag9 \\
	&=\ip{\lap_h\omega}{J_t^2\omega} - 2\sum_{i=1}^n \ips{J_t {\partial J_t\over \partial x^i}\grad_i \omega}{\omega} - \ip{R(\bar x)\omega}{J_t^2\omega} \pm c\|J_t\omega\|^2.\endalign$$
Now we estimate:
$$\align 2 \Big|\ips{J_t {\partial J_t\over \partial x^i}\grad_i \omega}{\omega}\Big| &\le  \sup_{\bar x} \left|{\partial J_t\over \partial x^i}\right|^2 \|J_t\grad_i\omega\|^2 +  \|\omega\|^2 \\
	&\le {C\over t^2} \|J_t\grad_i\omega\|^2 + C \|\omega\|^2\endalign$$
and
	$$\left|\ip{R(\bar x)\omega}{J_t^2\omega}\right| \le \sup_{\bar x} |R(\bar x)| \|J_t\omega\|^2 \le C \|\omega\|^2.$$
For large $t$, the two inequalities above, together with (9), give:
$$\align \sum_{i=1}^n \|J_t\grad_i\omega\|^2 + c^2\|x^1 J_t \omega\|^2 &\le C\left(\ip{\lap_h\omega}{J_t^2\omega} + \|J_t\omega\|^2\right) \tag10\\
	&\le C\left(\ip{\lap_h\omega}{\omega} + \|\omega\|^2\right).\endalign$$
G\"arding's inequality now follows if we take $t\to\infty$ in the inequality above. \qed\enddemo

\proclaim{21. Lemma}  Let $\omega$ be a smooth form. Suppose $\operatorname{supp}\omega \subset (1,\infty)\times U_j$ and suppose $\omega$, $\lap_h\omega\in L^2\Omega((1,\infty)\times Q)$.  Then
	$$\sum_{j=1}^n \|J_t \partial_j \omega\|^2 + \|J_t x^1 \omega\|^2 \le C\left(|\ip{J_t^2 \lap_h\omega}{\omega}| + \|J_t \omega\|^2\right),\tag11$$
and
	$$\sum_{j=1}^n \|\partial_j \omega\|^2 + \|x^1 \omega\|^2 \le C\left(\ip{\lap_h\omega}{\omega} + \|\omega\|^2\right).\tag12$$\endproclaim

\demo{Proof}  Recall that
	$$\grad_i(\omega_I dx^I) = (\grad_i \omega_I) dx^I + \omega_I \grad_i(dx^I) = (\partial_i \omega_I) dx^I + \omega_I \Gamma_i(\bar x)(dx^I), \tag13$$
where $\Gamma_i(\bar x):\wedge T_{\bar x}^* \to \wedge T_{\bar x}^*$ is an endomorphism.  The matrix of $\Gamma_i(\bar x)$ can be expressed in terms of the Christoffel symbols $\Gamma_{ij}^k(\bar x)$.

It follows from (13) that
	$$\|J_t\grad_j\omega\|^2 = \|J_t\partial_j\omega\|^2 + 2\ip{J_t\Gamma_j(\bar x)\omega}{J_t\partial_j\omega} + \|J_t\Gamma_j(\bar x)\omega\|^2.\tag14$$
For some choice of small $\eps>0$, we estimate:
	$$2|\ip{J_t\Gamma_j(\bar x)\omega}{J_t\partial_j\omega}| \le {1\over \eps} \|J_t\Gamma_j \omega\|^2 + \eps \|J_t\partial_j\omega\|^2 \le {C\over \eps} \|J_t\omega\|^2 + \eps \|J_t\partial_j\omega\|^2.$$

Therefore, from (14) 
	$$\|J_t\partial_j\omega\|^2 \le C \left(\|J_t\grad_j\omega\|^2 + \|J_t\omega\|^2\right).\tag15$$
Inequality (11) follows from (10) and (15).  Inequality (12) follows from (7) and (15) as $t \to \infty $.\qed\enddemo

\demo{Proof of theorem 16}  Induction argument.  For $\omega\in C_h^\infty$, $\operatorname{supp}\omega \subset (0,\infty)\times U_j$, inequality (4) follows from
	$$\|\omega\|_l^2 := \sum_{\matrix k,\alpha \\ k+|\alpha|\le l\endmatrix} \|(x^1)^k \partial^\alpha \omega\|^2 \quad \le \quad C(l) \left(\ip{\lap_h^l\omega}\omega + \|\omega\|^2\right)\tag16$$
where $0\le l<\infty$.

\proclaim{Remark} Operator $\lap_h^l$ is a PDO of order $2l$ with smooth coefficients which grow at most as 
$(x_1)^{2l}$. Therefore, for all $\omega \in C_h^{\infty }$ which satisfy $\|\omega \| \le \infty $,
we have 
$$\ip{\lap_h^l\omega}\omega \le C_l\|\omega \|^2_l. \tag 17$$
Moreover, for all $l_0\le l$
$$\ip{\lap_h^{l_0}\omega}\omega \le C\left( \ip{\lap_h^l\omega}\omega +|\ \omega \|^2\right).$$
\endproclaim

We will prove (16) by induction in $l$.  Lemma 21 proves (16) for $l=1$ and all $\omega\in C_h^\infty$ with $\operatorname{supp}\omega \subset (0,\infty)\times U_j$.

Suppose (16) is true for all $l\le l_0$ and all $\omega$ as above.  Let $k+|\alpha |= l_0$.  We substitute $(x^1)^k \partial^\alpha \omega$ instead of $\omega$ into (11) to get
$$\align \sum_{j=1}^n &\|J_t \partial_j (x^1)^k \partial^\alpha \omega\|^2 + \|J_t (x^1)^{k+1} \partial^\alpha \omega\|^2 \tag18\\
	&\le C\left(\left|\ip{J_t^2 \lap_h\left((x^1)^k \partial^\alpha\omega\right)}{(x^1)^k \partial^\alpha\omega}\right| + \|J_t (x^1)^k \partial^\alpha\omega\|^2\right).\endalign$$

Define $A=\ip{J_t^2 \lap_h\left((x^1)^k \partial^\alpha\omega\right)}{(x^1)^k \partial^\alpha\omega}$.  We write
	$$\lap_h\left((x^1)^k \partial^\alpha\omega\right) = (x^1)^k \partial^\alpha \lap_h\omega + \left[\lap_h,(x^1)^k \partial^\alpha\right]\omega,\tag19$$
where $\left[\cdot,\cdot\right]$ denotes the commutator of two differential operators.  Then from (6),
$$\align \left[\lap_h,(x^1)^k \partial^\alpha\right]\omega &= \left[\sum_{i=1}^n \grad_i^* \grad_i + R(\bar x),(x^1)^k \partial^\alpha\right]\omega + c^2\left[(x^1)^2,(x^1)^k \partial^\alpha\right]\omega \tag20\\
	&= P_1\omega + P_2\omega.\endalign$$

\proclaim{Proposition}

{\rm a)}  $P_1$ is a PDO of order at most $|\alpha|+1$ with coefficients which grow at most as $(x^1)^k$.

{\rm b)}  $P_2$ is a PDO of order at most $|\alpha|$ with coefficients which grow at most as $(x^1)^{k+1}$.\endproclaim

\demo{Proof of Proposition}  This follows from the formula (13) for $\grad_i$ and from two simple observations.

First, 
	$$\left[(x^1)^k,\partial_j\right] = \cases 0, &j\ne 1\\ k(x^1)^{k-1}, &j=1.\endcases$$

Second, the commutator of a smooth function $f(\bar x)$ with $\partial^\alpha$ is a PDO of order $|\alpha|-1$ with smooth (bounded) coefficients which depend only on $\bar x$. {\sl Q.E.D.}\enddemo

\smallskip
With the help of (19) and (20) we rewrite the quantity $A$ as
	$$A= \ip{J_t^2 (x^1)^k \partial^\alpha \lap_h\omega}{(x^1)^k \partial^\alpha\omega} + \ip{J_t^2 P_1 \omega}{(x^1)^k \partial^\alpha\omega} + \ip{J_t^2 P_2 \omega}{(x^1)^k \partial^\alpha\omega},\tag21$$
and we estimate:
$$\align |A| \le~&\|J_t (x^1)^k \partial^\alpha \lap_h\omega\| \cdot \|J_t (x^1)^k \partial^\alpha\omega\| \tag22 \\
	&+ \eps \|J_t P_1\omega\|^2 + {1\over\eps} \|J_t (x^1)^k \partial^\alpha\omega\|^2 + \eps \|J_t P_2\omega\|^2 + {1\over\eps} \|J_t (x^1)^k \partial^\alpha\omega\|^2,\endalign$$
where $\eps>0$ is to be chosen later.  Now we estimate the right--hand side of (22).  For $\lap_h\omega\in C_h^\infty$, by induction and by the remark above  we have
	$$\|J_t (x^1)^k \partial^\alpha \lap_h\omega\| \cdot \|J_t (x^1)^k \partial^\alpha\omega\| \le C\left(\ip{(\lap_h)^{l_0+1}\omega}\omega + \|\omega\|^2\right).$$
Moreover, from the proposition,
	$$\|J_t P_1\omega\|^2 + \|J_t P_2\omega\|^2 \le C \sum_{m+|\beta|\le l_0+1} \|J_t (x^1)^m \partial^\beta \omega\|^2.$$

Therefore,
	$$|A|\le C \Big(\ip{(\lap_h)^{l_0+1}\omega}\omega + \|\omega\|^2 + \eps\cdot \sum_{m+|\beta|\le l_0+1} \|J_t (x^1)^m \partial^\beta \omega\|^2 + \|\omega\|_{l_0}^2\Big).$$
The inequality (18) can be rewritten as
$$\align \sum_{j=1}^n &\|J_t (x^1)^k \partial_j \partial^\alpha \omega\|^2 + \|J_t (x^1)^{k+1} \partial^\alpha \omega\|^2 \tag23\\
	&\le C \Big(\ip{(\lap_h)^{l_0+1}\omega}\omega + \|\omega\|^2 + \|\omega\|_{l_0}^2\Big) + C \eps\cdot \sum_{m+|\beta|\le l_0+1} \|J_t (x^1)^m \partial^\beta \omega\|^2.\endalign$$
We can sum up inequalities (23) over all $k+|\alpha|=l_0$ and use the equality
$$\sum_{k+\alpha\le l_0} \sum_{j=1}^n \|J_t (x^1)^k \partial_j \partial^\alpha \omega\|^2 + \|J_t (x^1)^{k+1} \partial^\alpha \omega\|^2 = \sum_{m+|\beta|\le l_0+1} \|J_t (x^1)^m \partial^\beta \omega\|^2$$
to get (24) from (23):
$$\align \sum_{m+|\beta|\le l_0+1} \|J_t (x^1)^m \partial^\beta \omega\|^2 \le~&C(l_0+1) \Big(\ip{(\lap_h)^{l_0+1}\omega}\omega + \|\omega\|^2 + \|\omega\|_{l_0}^2\Big) \tag24 \\
	&+ \eps~C(l_0+1) \Big(\sum_{m+|\beta|\le l_0+1} \|J_t (x^1)^m \partial^\beta \omega\|^2\Big).\endalign$$
Finally, we choose $\eps$ small enough to move the right--hand $\sum_{m+|\beta|\le l_0+1}$ over to the left.

Then (16) follows from (24) for $l=l_0+1$ if we take $t\to\infty$. \qed \enddemo

\heading 6. Application: K\"unneth Formula \endheading

\proclaim{22. Proposition}  Let $M$, $N$ be orientable topologically tame manifolds.  Then
	$$H_c(M\times N) \cong H_c(M)\otimes H_c(N).$$\endproclaim

\demo{Proof}  Choose metrics $g_M$, $g_N$ so that $M$, $N$ have cylindrical ends as in theorem 2, and define $d\mu_M$, $d\mu_N$ as described there, choosing $c>0$.  Apply theorem 15 to $M$, $N$, and note that the same methods used to prove theorem 15 also show that $H_c(M\times N) \cong H_{\mu_M \times \mu_N}$, where $M\times N$ has been given the product metric and the product measure.

Now, $H_{\mu_M\times \mu_N}$ has the expected combinatorics, coming purely from analysis on Hilbert tensor products---more precisely Fredholm complexes.  See \BruningLesch.  Specifically, 
	$$H_{\mu_M\times \mu_N} \cong H_{\mu_M} \otimes H_{\mu_N},$$
which proves our theorem.\qed \enddemo

The above theorem applies only to what \BottTu~call manifolds with ``finite good cover'', but it is nonetheless a Hodge--theoretic proof of a K\"unneth formula for certain noncompact manifolds.  See \GriffithsHarris, chapter 0, for proof of the K\"unneth formula on compact manifolds, using the same Hodge--theoretic technique of proof.  The K\"unneth formula in fact holds for any pair of manifolds $M$, $F$ for which the cohomology of $F$ is finite--dimensional.

\enddocument